\documentclass[oneside,english]{amsart}
\usepackage[T1]{fontenc}
\usepackage[latin9]{inputenc}
\usepackage{amsthm}
\usepackage{amstext}
\usepackage{amssymb}
\usepackage{esint}

\makeatletter
\numberwithin{equation}{section}
\numberwithin{figure}{section}
\theoremstyle{plain}
\newtheorem{thm}{\protect\theoremname}[section]

\theoremstyle{plain}
\theoremstyle{definition}

\theoremstyle{plain}

\theoremstyle{plain}
\newtheorem{rem}[thm]{\protect\remarkname}
\theoremstyle{plain}
\makeatother

\usepackage{babel}
\usepackage{color}
\providecommand{\definitionname}{Definition}
\providecommand{\lemmaname}{Lemma}
\providecommand{\theoremname}{Theorem}
\providecommand{\corollaryname}{Corollary}
\providecommand{\remarkname}{Remark}
\providecommand{\propositionname}{Proposition}

\DeclareMathOperator{\loc}{loc}

\DeclareMathOperator{\ess}{ess}
\DeclareMathOperator{\cp}{cap}

\begin{document}

\title[On the boundary behavior of weak $(p,q)$-quasiconformal mappings]
{On the boundary behavior of weak $(p,q)$-quasiconformal mappings}

\author{Vladimir Gol'dshtein, Evgeny Sevost'yanov, Alexander Ukhlov}
\begin{abstract}
Let $\Omega$ and $\widetilde{\Omega}$ be domains in the Euclidean space $\mathbb R^n$. We study the boundary behavior of weak $(p,q)$-quasiconformal mappings, $\varphi:\Omega\to \widetilde{\Omega}$, $n-1<q\leq p<n$. The suggested method is based on the capacitary distortion properties of the weak $(p,q)$-quasiconformal mappings.
\end{abstract}
\maketitle

\section{Introduction }

Let $\Omega$ and $\widetilde{\Omega}$ be domains in the Euclidean space $\mathbb R^n$, $n\geq 2$. We consider the boundary behavior of weak $(p,q)$-quasiconformal mappings $\varphi:\Omega\to \widetilde{\Omega}$, $n-1<q\leq p<n$. The weak $(p,q)$-quasiconformal mappings represent generalizations of (quasi)conformal mappings and have significant applications in the geometric analysis of PDE (see, for example, \cite{GPU18_3,GU09,GU16,GU17}). In the frameworks of the boundary value problems for elliptic equations becomes important the boundary behavior of this class of mappings. From the historic point of view it arises to the boundary behavior of conformal mappings (univalent analytic functions) \cite{C13} and to the boundary behavior of quasiconformal mappings \cite{Ahl66}. In series of works \cite{GV78,KR16,KP87,Na70,O66,Smol10,S85,Z62} the boundary behavior of space quasiconformal mappings and their generalizations in the terms of  capacitory (moduli) inequalities was studied.

The theory of weak $(p,q)$-quasiconformal mappings arose in the Sobolev embedding theory \cite{GS82,M69} and was founded in the series of works \cite{GGR95,U93,VU98,VU02,VU04}. Recent applications of the weak $(p,q)$-quasiconformal mappings theory to Sobolev extension operators can be found in \cite{KZ22,KUZ22}. Recall that a mapping $\varphi:\Omega\to\widetilde{\Omega}$ is called a weak $(p,q)$-quasiconformal mapping if $\varphi\in W^1_{1,\loc}(\Omega)$, has finite distortion and
$$
K_{p,q}^{\frac{pq}{p-q}}(\varphi;\Omega)=\int\limits_{\Omega}\left(\frac{|D\varphi(x)|^p}{|J(x,\varphi)|}\right)^{\frac{q}{p-q}}~dx<\infty,
$$
for $1< q<p< \infty$ \cite{U93,VU98}  and
$$
K_{p,p}^{p}(\varphi;\Omega)=\ess\sup\limits_{\Omega}\frac{|D\varphi(x)|^p}{|J(x,\varphi)|}<\infty,
$$
for $1< q=p< \infty$ \cite{GGR95,V88}. In the case $p=q=n$ we have quasiconformal mappings \cite{V71} and in the case $1<q<p=n$ we have mappings of integrable distortion \cite{HK93}.

The main result of the article states: {\it Let $\varphi :\Omega\to
\widetilde{\Omega},$ $\varphi(\Omega)=\widetilde{\Omega},$ be a weak
$(p,q)$-quasiconformal mapping, $n-1<q\leq p<n,$ where $\Omega$ has
a strongly accessible boundary with respect to $q$-capacity and
$\widetilde{\Omega}$ has locally connected boundary. Then the
inverse mapping $\varphi^{-1}$ can be extended to a continuous
mapping
$$
\overline{\varphi^{-1}}: \overline{\widetilde{\Omega}}\to \overline{\Omega}.
$$}
The definition of a strongly accessible boundary in the terms of the $q$-capacity will be given in Section 3.

The suggested method is based on the capacitary distortion
properties of the weak $(p,q)$-quasiconformal mappings. Let us note that the weak $(p,q)$-quasi\-con\-for\-mal mappings are closely connected with mappings defined by
weighted capacitary (moduli) inequalities \cite{GSU22,MU21_2}.

\section{Weak quasiconformal mappings}

Let $\Omega$ be a domain in the Euclidean space $\mathbb R^n$, $n\geq 2$. The Sobolev space $W^1_p(\Omega)$, $1\leq p\leq\infty$,
is defined
as a Banach space of locally integrable weakly differentiable functions
$f:\Omega\to\mathbb{R}$ equipped with the following norm:
\[
\|f\mid W^1_p(\Omega)\|=\| f\mid L_p(\Omega)\|+\|\nabla f\mid L_p(\Omega)\|,
\]
where $\nabla f$ is the weak gradient of the function $f$.

The homogeneous seminormed Sobolev space $L^1_p(\Omega)$, $1\leq p\leq\infty$, is defined as a space
of locally integrable weakly differentiable functions $f:\Omega\to\mathbb{R}$ equipped
with the following seminorm:
\[
\|f\mid L^1_p(\Omega)\|=\|\nabla f\mid L_p(\Omega)\|.
\]

In accordance with the non-linear potential theory \cite{MH72} we consider elements of Sobolev spaces $W^1_p(\Omega)$ are equivalence classes up to a set of $p$-capacity zero  \cite{M}. 

Let $\Omega$ and $\widetilde{\Omega}$ be domains in $\mathbb R^n$, $n\geq 2$. We say that
a homeomorphism $\varphi:\Omega\to\widetilde{\Omega}$ induces a bounded composition
operator
\[
\varphi^{\ast}:L^1_p(\widetilde{\Omega})\to L^1_q(\Omega),\,\,\,1\leq q\leq p\leq\infty,
\]
by the composition rule $\varphi^{\ast}(f)=f\circ\varphi$, if for
any function $f\in L^1_p(\widetilde{\Omega})$, the composition
$\varphi^{\ast}(f)\in L^1_q(\Omega)$ is defined quasi-everywhere in
$\Omega$ \cite{M} and there exists a constant
$K_{p,q}(\Omega)<\infty$ such that
\[
\|\varphi^{\ast}(f)\mid L^1_q(\Omega)\|\leq K_{p,q}(\Omega)\|f\mid L^1_p(\widetilde{\Omega})\|.
\]

In the geometric function theory composition operators on Sobolev spaces arise in the work \cite{VG75} and have numerous applications in the geometric analysis of PDE. The characterization of composition operators on Sobolev spaces
is given in the following theorem \cite{U93,VU98} (\cite{V88} for the case $p=q>n$ and \cite{GGR95} for the case $n-1<q=p<n$).

\begin{thm}
Let $\varphi:\Omega\to\widetilde{\Omega}$ be a homeomorphism
between two domains $\Omega$ and $\widetilde{\Omega}$. The homeomorphism  $\varphi$ induces a bounded composition
operator
\[
\varphi^{\ast}:L^1_p(\widetilde{\Omega})\to L^1_q(\Omega),\,\,\,1\leq q<p<\infty,
\]
if and only if $\varphi$ is the weak $(p,q)$-quasiconformal mapping.
\end{thm}

Recall the notion of a variational $p$-capacity \cite{VGR79}. The condenser in a domain $\Omega\subset \mathbb R^n$ is a pair $(E,F)$ of connected disjoint closed relatively to $\Omega$ sets $E,F\subset \Omega$. A continuous function $f\in L_p^1(\Omega)$ is called an admissible function for the condenser $(E,F)$,
if the set $E\cap \Omega$ is contained in some connected component of the set $\operatorname{Int}\{x\vert u(x)=0\}$ and the set $F\cap \Omega$ is contained in some connected component of the set $\operatorname{Int}\{x\vert u(x)=1\}$. We call as the $p$-capacity of the condenser $(E,F)$ relatively to domain $\Omega$
the following quantity:
$$
{\cp}_p(E,F;\Omega)=\inf\|u\vert L_p^1(\Omega)\|^p.
$$
Here the greatest lower bond is taken over all functions admissible for the condenser $(E,F)\subset\Omega$. If the condenser has no admissible functions we put the capacity equal to infinity.

The following capacitory properties of weak $(p,q)$-quasiconformal mappings were established in \cite{U93,VU98} (for the case $p=\infty$ see in \cite{U04}).

\begin{thm}
\label{theorem:CapacityDescPQ_O}
Let $\varphi :\Omega\to \widetilde{\Omega}$ be a weak $(p,q)$-quasiconformal mapping, $1\leq q\leq p\leq \infty$. Then
for every condenser
$(E,F)\subset \widetilde{\Omega}$
the inequality
$$
\cp_{q}^{1/q}(\varphi^{-1}(E),\varphi^{-1}(F);\Omega)
\leq K_{p,q}(\varphi;\Omega) \cp_{p}^{1/p}(E,F;\widetilde{\Omega})
$$
holds.
\end{thm}

By using these capacitory distortion properties we consider boundary behavior of weak $(p,q)$-quasiconformal mappings.

\section{On boundaries correspondence }

Recall that the boundary $\partial \Omega$ of a~domain $\Omega\subset\mathbb R^n$ is called {\it 
strongly accessible at a~point $x_0\in \partial \Omega$ with respect
to the $p$-capacity} if for each neighborhood $U$ such that $\partial U \cap \Omega \neq \emptyset$ of $x_0$ there
exist a~compact set $E\subset \Omega$, a~neighborhood $V\subset U$ of
$x_0$ and $\delta>0$ such that
\begin{equation}
\label{eq1.3_a} \cp_p(E, F; \Omega)\geqslant \delta
\end{equation}
for each continuum $F$ in~$\Omega$ that intersects $\partial U$ and
$\partial V$. The boundary $\partial \Omega$ of a~domain $\Omega$ is called
{\it strongly accessible at a~point $x_0\in \partial \Omega$}, if it is
strongly accessible at a~point $x_0\in \partial \Omega$ with respect to
the $n$-capacity ($n$-modulus).

\begin{rem}\label{rem1}
The notion of a strongly accessible boundary was introduced
in~\cite[section~3.8]{MRSY09} and it is very close to the notion of a
uniform domain which was introduced by N\"{a}kki~\cite{Na73} Theorem 6.2.
\end{rem}

This notion also coincides, up to some (not too essential) details, with the concept of a quasiconformally accessible
boundary~\cite[section~1.7]{Na70}. Note that both of these concepts were formulated in terms of the modulus of families of paths. In
this connection, we recall the concept of a modulus of a family of locally rectifiable curves (paths) $\Gamma $.

\medskip
Let $\rho:{\mathbb R}^n \rightarrow [0,\infty]$ be a Borel function. Then $\rho$  is called
{\it admissible} for $\Gamma$ (i.~e. $\rho \in {\rm adm}\, \Gamma)$, if the inequality $\int\limits_{\gamma} \rho(x)~ds \geqslant 1$ holds for any locally rectifiable curve $\gamma\in\Gamma$. Let $p\geqslant 1$, the quantity
\begin{equation}\label{eq2_a}
M_p(\Gamma)\,=\,\inf\limits_{\rho \in {\rm adm}\, \Gamma}
\int\limits_{{\mathbb R}^n} \rho^p(x)\, dx
\end{equation}
is called the $p$-modulus of the family of curves of $\Gamma$.

For a given domain $\Omega$ in $\overline{{\mathbb R}^n}={\mathbb R}^n\cup \{\infty\},$ and sets  $E$ and $F$ in $\Omega$
we denote by the symbol $\Gamma(E, F, \Omega)$ the family of all locally rectifiable curves $\gamma:[0, 1]\rightarrow
\overline{{\mathbb R}^n}$ joining the sets $E$ and $F$ in $\Omega$, that is: $\gamma(0)\in E,$ $\gamma(1)\in F$ and $\gamma(t)\in \Omega$ for all $t\in (0, 1)$.

We say that the boundary $\partial \Omega$ of a domain $\Omega$ is {\it
strongly accessible at a~point $x_0\in \partial \Omega$ with respect to
the $p$-modulus} if for each neighborhood $U$ of $x_0$ there exist
a~compact set $E\subset \Omega$, a neighborhood $V\subset U$ of $x_0$ and
$\delta>0$ such that
\begin{equation}\label{eq1A} M_p(\Gamma(E, F, \Omega)\geqslant \delta
\end{equation}
for each continuum $F$ in $\Omega$ that intersects $\partial U$ and
$\partial V$.

Note that (\ref{eq1.3_a}) is equivalent to (\ref{eq1A}). Indeed,
$$M_p(\Gamma(E, F, \Omega)=\cp_p(E, F;\Omega)$$
by Hesse and Shlyk equalities (see \cite[Theorem~5.5]{Hes75} and~\cite[Theorem~1]{Shl93}).

\medskip
Based on the definition of domains with strongly accessible boundaries and taking Remark~\ref{rem1} into account, we give
some examples of such domains.

\medskip
\noindent
1. By Theorem~6.2 and Corollary~6.8 in \cite{Na73}, the planar domain with finitely many boundary components has a strongly
accessible boundary whenever it is finitely connected on the boundary.

\medskip
\noindent
2. Following \cite{Na73}, a domain $\Omega$ is said to be {\it quasiconformally collared} on the boundary if each point of
$\partial \Omega$ has arbitrarily small neighborhoods $U$ such that $U\cap \Omega$ can be mapped quasiconformally onto a ball $B\subset\mathbb R^n$. Let $\Omega$ be a domain which can be mapped quasiconformally onto some quasiconformally collared domain. If $\Omega$ is finitely connected on the boundary, then the boundary of $\Omega$ is strongly accessible (see Corollary~6.7 in \cite{Na73}).

\medskip
\noindent
3. The next example gives domains with a strongly accessible boundary for $p\ne n$.
Recall the notion of the upper gradient \cite{HK00,He01}. Let $(X, \mu)$ be a metric measure space. A Borel function
$g\colon X\rightarrow [0, \infty]$ is said to be an {\it upper
gradient} of a function $u:X\to{\mathbb R}$ if
$$
|u(x)-u(y)|\leqslant
\int\limits_{\gamma}g\,ds
$$
for any rectifiable curve $\gamma$ joining $x$ and $y$ in $X$.
Let $1\leqslant p<\infty.$ We say that $X, \mu$ admits {\it a $(1;
p)$-Poincare inequality} if there exists a constant $1\leq C_p<\infty$ such
that
$$
\frac{1}{\mu(B)}\int\limits_{B}|u-u_B|d\mu(x)\leqslant C_p\cdot({\rm
diam\,}B)\left(\frac{1}{\mu(B)} \int\limits_{B}g^p
d\mu(x)\right)^{1/p}
$$
for all balls $B$ in $X$ and for all bounded continuous functions $u$ on $B$, where $g$ is the upper gradient of $u$. Metric measure
spaces where the inequalities
$$
\frac{1}{K}R^{n}\leqslant \mu(B(x_0,R))\leqslant KR^{n}
$$
hold with some constant $1\leq K<\infty$, every $x_0\in X$  and all
$R<{\rm diam}\,X$ are called {\it Ahlfors $n$-regular.}
Let $\Omega\subset B(0, R)$ for some $R>0$ be an Ahlfors $n$-regular domain satisfying the $(1;p)$-Poincare inequality for $n-1<p\leqslant n$.
Assume that $E$ and $F$ are some continua $E,F\subset \Omega$. By \cite[Proposition~4.7]{AS10} and due to Remark~\ref{rem1},
\begin{equation}\label{eq1H}
\cp_p(E, F; \Omega)\geqslant  \frac{1}{C}\frac{\min\{{\rm diam} E, {\rm
diam} F\}}{R^{1+p-n}}\,,
\end{equation}
where $C>0$ is some constant.

\medskip
\noindent
Let us prove that $\Omega$ has a strongly accessible boundary at any point $x_0\in \partial \Omega$ with respect to $p$-capacity. Suppose that $x_0\in \partial \Omega$ and that $U$ is an arbitrary neighborhood of $x_0$. Choose a small enough $\varepsilon_1>0$ such that for $V:=B(x_0, \varepsilon_1)$, we have $\overline {V}\subset \overline{U}$.  Because $(\partial U \cap D) \ne\varnothing$ we can set $\varepsilon_2:=dist(\partial U,\partial V)>0$. Note that for arbitrary continua $F_1$ and $F_2$ in
$\Omega$ satisfying $F_1\cap \partial U \ne\varnothing\ne F_1\cap\partial V$ and $F_2\cap \partial U \ne\varnothing\ne F_2\cap\partial V$ we
have ${\rm diam} (F_1)\geqslant \varepsilon_2$ and ${\rm diam}(F_2)\geqslant \varepsilon_2$. Therefore, by~(\ref{eq1H}), we
obtain $\cp_p(\Gamma(F_1, F_2; G_0))\geqslant \varepsilon_2$, as required.

\medskip
Given $x_0\in {\mathbb R}^n,$ we set
$$
B(x_0, r)=\{x\in {\mathbb R}^n: |x-x_0|<r\}\,,$$
\begin{equation}\label{eq1ED} S(x_0,r) = \{ x\,\in\,{\mathbb R}^n : |x-x_0|=r\}\,,
\end{equation}
\begin{equation}\label{eq1**A} A=A(x_0, r_1, r_2)=\{ x\,\in\,{\mathbb R}^n :
r_1<|x-x_0|<r_2\}\,.
\end{equation}

\medskip
The following statement holds.

\begin{thm}
\label{theorem:Bound_PQ_G} Let $\varphi :\Omega\to \widetilde{\Omega},$ $\varphi(\Omega)=\widetilde{\Omega}$, be a weak
$(p,q)$-quasiconformal mapping, $n-1<q\leq p<n,$ $\Omega$ has a strongly accessible boundary with a respect to $q$-capacity and
$\widetilde{\Omega}$ has locally connected boundary. Then the inverse mapping $\varphi^{-1}$ can be extended by continuity to
the continuous mapping
$$
\overline{\varphi^{-1}}: \overline{\widetilde{\Omega}}\to \overline{\Omega}.
$$
\end{thm}

\medskip
\begin{proof}
Suppose the contrary, namely, there exists such point $b\in \partial \widetilde{\Omega}$, that $\varphi^{-1}$ has no a continuous
extension to the point $b$. It means that there exist two sequences $x_i$, $x_i^{\,\prime}\in \widetilde{\Omega}$, $i=1,2,\ldots$, such that
$x_i\rightarrow b,$ $x_i^{\,\prime}\rightarrow b$ as $i\rightarrow
\infty,$ and $\varphi^{\,-1}(x_i)\rightarrow y,$
$\varphi^{\,-1}(x_i^{\,\prime})\rightarrow y^{\,\prime}$ as
$i\rightarrow \infty,$ while $y^{\,\prime}\ne y.$ Note that $y$ and
$y^{\,\prime}\in \partial \Omega$, because $C(\varphi^{-1},
\partial \widetilde{\Omega})\subset \partial \Omega$ for any homeomorphism $\varphi^{-1}$,
see \cite[Proposition~13.5]{MRSY09}.

Here
$$
C(\varphi^{-1}, \partial \widetilde{\Omega})=\bigcup\limits_{x_0\in \partial\widetilde{\Omega}} C(\varphi^{-1}, x_0)\,,
$$
where
$$
C(\varphi^{-1}, x_0)=\{y\in\overline{{\mathbb R}^n}:
\,\,\exists\,\,x_k\in \widetilde{\Omega}, x_k\rightarrow x_0:
\varphi^{\,-1}(x_k)\rightarrow y, k\rightarrow\infty\}.
$$

By the definition of a strongly accessible boundary at the point $y\in \partial \Omega$ with respect to the $q$-capacity, for any
neighborhood $U$ of this point there exists a compact set $C_0^{\,\prime}\subset \Omega$, a neighborhood $V$ of a point $y$,
$V\subset U$, and a number $\delta>0$ such that
\begin{equation}\label{eq1}
\cp_q(C_0^{\,\prime}, F; \Omega)\geqslant \delta
>0
\end{equation}
for any continua $F,$ intersecting $\partial U$ and $\partial V.$
Since $C(\varphi^{\,-1},
\partial \widetilde{\Omega})\subset \partial \Omega,$ we obtain that the condition $C_0\cap \partial
\widetilde{\Omega}=\varnothing$ holds for
$C_0:=\varphi(C_0^{\,\prime})$. Now suppose $\varepsilon_0>0$ is such that
$C_0\cap\overline{B(b, \varepsilon_0)}=\varnothing$.

Since $\widetilde{\Omega}$ is locally connected at $b,$ we may join
the points $x_i$ and $x_i^{\,\prime}$ by a path $\gamma_i$ lying in
$V\cap \widetilde{\Omega}.$ We may consider that $\gamma_i\in
\overline{B(b, 2^{\,-i})}\cap \widetilde{\Omega}.$ Since
$\varphi^{\,-1}(x_i)\in V$ and $\varphi^{\,-1}(x_i^{\,\prime})\in
\widetilde{\Omega}\setminus \overline{U}$ for sufficiently large
$i\in {\mathbb N}$ and due to~(\ref{eq1}), we may find $i_0\in {\mathbb
N}$ such that 
\begin{equation}\label{eq2}
\cp_q(C_0^{\,\prime}, \varphi^{\,-1}(\gamma_i);
\Omega))\geqslant \delta
>0
\end{equation}
for any $i\geqslant i_0\in {\mathbb N}.$ Immerse the compact $C_0$ into
some the continuum $C_1,$ still completely lying in the domain
$\widetilde{\Omega},$ see~\cite[Lemma~1]{Smol10}. By reducing
$\varepsilon_0>0,$ we may again assume that $C_1\cap\overline{B(b,
\varepsilon_0)}=\varnothing.$
By Theorem~\ref{theorem:CapacityDescPQ_O} 
\begin{multline}
\label{eq7}
\cp^{1/q}_q C_0^{\,\prime}, \varphi^{\,-1}(\gamma_i);
\Omega)\leqslant \\
\cp_{q}^{1/q}(\varphi^{-1}(\gamma_i),
\varphi^{\,-1}(C_1);\Omega) \leqslant K_{p,q}(\varphi;\Omega)
\cp_{p}^{1/p}(\gamma_i, C_1;\widetilde{\Omega})\,.
\end{multline}

\medskip
Let us prove that $\cp_{p}(\gamma_i,
C_1;\widetilde{\Omega})\rightarrow 0$ as $i\rightarrow\infty.$
Indeed, by the definition of the capacity
\begin{equation}\label{eq8}
\cp_p(\gamma_i, C_1;\widetilde{\Omega})\leqslant
\cp_p(S(b, 2^{\,-i}), S(b, \varepsilon_0); A(b, 2^{\,-i},
\varepsilon_0))\,.
\end{equation}
However, by \cite[relation~(2)]{Geh71} and due to the Remark~\ref{rem1}
$$
\cp_n(S(b, 2^{\,-i}), S(b, \varepsilon_0); A(b, 2^{\,-i},
\varepsilon_0))=\frac{\omega_{n-1}}{\log^{n-1}\frac{\varepsilon_0}{2^{\,-i}}}\rightarrow
0\,,\qquad i\rightarrow\infty\,,
$$
and
\begin{multline*}
\cp_p(S(b, 2^{\,-i}), S(b, \varepsilon_0); A(b, 2^{\,-i},
\varepsilon_0))\\=
\left(\frac{n-p}{p-1}\right)^{p-1}\frac{\omega_{n-1}}
{\left((2^{\,-i})^\frac{p-n}{p-1}-\varepsilon_0^\frac{p-n}{p-1}\right)^{p-1}}\rightarrow
0\,,\quad i\rightarrow \infty \quad p\ne n\,.
\end{multline*}
Now, $\cp_{p}(\gamma_i, C_1;\widetilde{\Omega})\rightarrow
0$ as $i\rightarrow\infty,$ as required. In this case, the
relation~(\ref{eq7}) contradicts with~(\ref{eq2}). This
contradiction proves the theorem.
\end{proof}

\medskip
In the case $n<q\leq p<\infty$ we use the following composition duality theorem \cite{U93,VU98}:

\begin{thm}
\label{CompThD} Let $\varphi :\Omega\to \widetilde{\Omega}$ be a weak $(p,q)$-quasiconformal mapping, $n-1<q\leq p<\infty$.
Then the inverse mapping $\varphi^{-1}:\widetilde{\Omega}\to\Omega$ generates a bounded composition operator 
\[
\left(\varphi^{-1}\right)^{\ast}:L^1_{q'}(\Omega)\to L^1_{p'}(\widetilde{\Omega}),
\]
where $p'=p/(p-n+1)$, $q'=q/(q-n+1)$. 
\end{thm}

By using this composition duality theorem we obtain the capacitory distortion estimate:

\begin{thm}
\label{theorem:CapacityDescPQ_O_I}
Let $\varphi :\Omega\to \widetilde{\Omega}$ be a weak $(p,q)$-quasiconformal mapping, $n<q\leq p<(n-1)^2/(n-2)$. Then
for every condenser
$(F_0,F_1)\subset {\Omega}$
the inequality
$$
\cp_{p'}^{1/p'}(\varphi(F_0),\varphi(F_1);\widetilde{\Omega})
\leq K_{q',p'}(\varphi^{-1};\widetilde{\Omega}) \cp_{q'}^{1/q'}(F_0,F_1;{\Omega}),\,\,\,\,n-1<p'\leq q'<n,
$$
holds, where $p'=p/(p-(n-1))$ and $q'=q/(q-(n-1))$.
\end{thm}

The condition $p<<(n-1)^2/(n-2)$ provides that $p'>n-1$. Hence, by using Theorem~\ref{theorem:Bound_PQ_G} and Theorem~\ref{theorem:CapacityDescPQ_O_I} we obtain:

\begin{thm}
\label{theorem:Bound_PQ_D}
Let $\varphi :\Omega\to \widetilde{\Omega}$ be a weak $(p,q)$-quasiconformal mapping, $n<q\leq p<(n-1)^2/(n-2)$.
Suppose that $\Omega$ has a locally connected boundary and
$\widetilde{\Omega}$ has strongly accessible boundary with a respect to $p'$-capacity, $p'=p/(p-n+1)$. Then the mapping $\varphi$ can be extended by continuity to
the continuous mapping
$$
\overline{\varphi}: \overline{\Omega}\to\overline{\widetilde{\Omega}}.
$$
\end{thm}

\vskip 0.5cm

Vladimir Gol'dshtein; Department of Mathematics, Ben-Gurion University of the Negev, P.O.Box 653, Beer Sheva, 8410501, Israel 
 
\emph{E-mail address:} \email{vladimir@math.bgu.ac.il} \\     

Evgeny Sevost'yanov; Department of Mathematical Analysis, Zhytomyr Ivan Franko State University, 40 Bol'shaya Berdichevskaya Str., Zhytomyr, 10008, Ukraine   

Institute of Applied Mathematics and Mechanics of NAS of Ukraine, 1 Dobrovol'skogo Str., 84 100 Slavyansk, Ukraine   
 
\emph{E-mail address:} \email{esevostyanov2009@gmail.com} \\    

Alexander Ukhlov; Department of Mathematics, Ben-Gurion University of the Negev, P.O.Box 653, Beer Sheva, 8410501, Israel 
							
\emph{E-mail address:} \email{ukhlov@math.bgu.ac.il


\medskip

\begin{thebibliography}{References}

\bibitem{AS10} T.~Adamowicz, N.~Shanmugalingam, Non-conformal Loewner type estimates for modulus of curve families,
Ann. Acad. Sci. Fenn. Math., 35 (2010), 609--626.

\bibitem{Ahl66}
L.~Ahlfors, Lectures on quasiconformal mappings. D. Van Nostrand Co., Inc., Toronto, Ont.-New York-London, 1966.

\bibitem{C13} C.~Carath\'eodory, \"Uber die Begrenzung einfach zusammenh\"angender Gebiete,
Math. Ann., 73 (1913), 323--370.

\bibitem{Geh71}
F.W.~Gehring, Lipschitz mappings and the $p$-capacity of rings in $n$-space, Advances in the Theory of Riemann Surfaces,
Princeton, University Press, 175--193, 1971.

\bibitem{GGR95} V.~Gol'dshtein, L.~Gurov, A.~Romanov, Homeomorphisms that induce monomorphisms of Sobolev spaces,
Israel J. Math., 91 (1995), 31--60.

\bibitem{GPU18_3} V.~Gol'dshtein, V.~Pchelintsev, A.~Ukhlov, On the First Eigenvalue of the Degenerate p-Laplace Operator in Non-convex Domains. Integral Equations Operator Theory, 90 (2018), 90:43.

\bibitem{GSU22} V.~Gol'dshtein, E.~Sevost'yanov, A.~Ukhlov, Composition operators on Sobolev spaces and weighted moduli inequalities, Math. Reports (accepted).

\bibitem{GS82} V.~M.~Gol'dshtein, V.~N.~Sitnikov, Continuation of functions of the class $W^1_p$ across H\"older boundaries, Imbedding theorems and their applica-ions, Trudy Sem. S. L. Soboleva, 1 (1982), 31--43.

\bibitem{GU09} V.~Gol'dshtein, A.~Ukhlov, Weighted Sobolev spaces and embedding theorems, Trans. Amer. Math. Soc., 361 (2009), 3829--3850.

\bibitem{GU16} V.~Gol'dshtein, A.~Ukhlov, On the first Eigenvalues of Free Vibrating Membranes in Conformal Regular Domains,
Arch. Rational Mech. Anal., 221 (2016), no. 2, 893--915.

\bibitem{GU17} V.~Gol'dshtein, A.~Ukhlov, The spectral estimates for the Neumann-Laplace operator in space domains,
Adv. in Math., 315 (2017), 166--193.

\bibitem{GV78} V.~M.~Gol'dshtein, S.~K.~Vodop'yanov, Metric completion of a domain by means of a conformal capacity that is invariant under quasiconformal mappings, Dokl. Akad. Nauk SSSR, 238 (1978), 1040--1042.

\bibitem{HK00} P.~Hajlasz, P.~Koskela, Sobolev met Poincar\'e, Mem. Amer. Math. Soc., 145 (2000), 1--101. 

\bibitem{He01} J.~Heinonen, Lectures on analysis on metric spaces, Springer Science+Business Media, New York, 2001.

\bibitem{HK93} J.~Heinonen, P.~Koskela, Mappings with integrable dilatation, Arch. Rational Mech. Anal., 125 (1993), 81--97.

\bibitem{Hes75} J.~Hesse, A $p$-extremal length and $p$-capacity equality, Ark. Mat., 13 (1975), 131--144.

\bibitem{KZ22} P.~Koskela, Z.~Zhu, Sobolev Extensions Via Reflections, J. Math. Sci., 268 (2022), 376--401.

\bibitem{KUZ22} P.~Koskela, A.~Ukhlov, Z.~Zhu, The volume of the boundary of a Sobolev $(p,q)$-extension domain, J. Funct. Anal., 283 (2022), 109703.

\bibitem{KR16} D.~Kovtonyuk, V.~Ryazanov, On boundary behavior of spatial mappings, Rev. Roumaine Math. Pures Appl., 61 (2016), 57--73.

\bibitem{KP87} V.~I.~Kruglikov, V.~I.~Paikov, Capacity and prime ends of a space domain, Dokl.~Akad.~Nauk Ukr.~SSR, Ser. A 5, 154 (1987),
10--13.

\bibitem{MRSY09} O.~Martio, V.~Ryazanov, U.~Srebro and E.~Yakubov, Moduli in Modern Mapping Theory, Springer Science + Business Media,
LLC, New York, 2009.

\bibitem{M69} V.~G.~Maz'ya, Weak solutions of the Dirichlet and Neumann problems,
Trudy Moskov. Mat. Ob-va., 20 (1969), 137--172 (1969)

\bibitem{M} V.~Maz'ya, Sobolev spaces: with applications to elliptic partial differential equations, Springer, Berlin/Heidelberg, 2010.

\bibitem{MH72} V.~G.~Maz'ya, V.~P.~Havin, Non-linear potential theory, Russian Math. Surveys, \textbf{27} (1972), 71--148.

\bibitem{MU21_2} A.~Menovschikov, A.~Ukhlov, Composition operators on Sobolev spaces, capacity and weighted Sobolev inequalities, arXiv:2110.09261.

\bibitem{Na70} R.~N\"{a}kki, Boundary behavior of quasiconformal mappings in $n$-space, Ann. Acad. Sci. Fenn. Ser. A., 484 (1970), 1--50.

\bibitem{Na73} R.~N\"{a}kki, Extension of Loewner's capacity theorem, Trans. Amer. Math. Soc., 180 (1973) , 229--236.

\bibitem{O66} I.~S.~Ovchinnikov, Prime ends of a certain class of space regions, Trudy Tomsk. Gos. Univ. Ser. Meh.-Mat., 189 (1966), 96--103.

\bibitem{Shl93} V.A.~Shlyk, The equality between $p$-capacity and $p$-modulus, Siberian Math. J., 34 (1993), 1196--1200.

\bibitem{Smol10} E.S.~Smolovaya, Boundary behavior of ring $Q$-homeomorphisms in metric spaces, Ukr. Math. J., 62 (2010), no.~5, 785--793.

\bibitem{S85} G.~D.~Suvorov, The generalized "length and area principle" in mapping theory, Naukova Dumka, Kiev, 1985.

\bibitem{U93} A.~Ukhlov, On mappings, which induce embeddings of Sobolev spaces, Siberian Math. J., 34 (1993), 185--192.

\bibitem{U04} A.~Ukhlov, Differential and geometrical properties of Sobolev mappings, Matem. Notes., 75 (2004), 291--294.

\bibitem{V71} J.~V\"ais\"al\"a, Lectures on $n$-dimensional quasiconformal mappings. Lecture Notes in Math. 229, Springer Verlag, Berlin, 1971.

\bibitem{V88} S.~K.~Vodop'yanov, Taylor Formula and Function Spaces, Novosibirsk Univ. Press., 1988.

\bibitem{VG75} S.~K.~Vodop'yanov, V.~M.~Gol'dshtein, Structure isomorphisms of spaces $W^1_n$ and quasiconformal mappings,
Siberian Math. J., 16 (1975), 224--246.

\bibitem{VGR79} S.~K.~Vodop'yanov, V.~M.~Gol'dshtein, Yu.~G.~Reshetnyak, On geometric properties of functions with generalized first derivatives, Uspekhi Mat. Nauk, 34,  17--65 (1979)

\bibitem{VU98} S.~K.~Vodop'yanov, A.~D.~Ukhlov, Sobolev spaces and $(P,Q)$-quasiconformal mappings of Carnot groups,
Siberian Math. J., 39 (1998), 665--682.

\bibitem{VU02} S.~K.~Vodop'yanov, A.~D.~Ukhlov, Superposition operators in Sobolev spaces, Russian Mathematics (Izvestiya VUZ), 46 (2002), no. 4, 11--33.

\bibitem{VU04} S.~K.~Vodop'yanov, A.~D.~Ukhlov, Set Functions and Their Applications in the Theory of Lebesgue and Sobolev Spaces. I, Siberian Adv. Math., 14:4 (2004), 78--125.

\bibitem{Z62} V.~A.~Zorich,  On boundary correspondence for $Q$-quasiconformal mappings of a sphere, Dokl. Akad. Nauk SSSR, 145 (1962), 31--34.




\end{thebibliography}
\end{document}